\documentclass[12pt]{article}
\newcommand{\be}{\begin{equation}}
\newcommand{\ee}{\end{equation}}
\newcommand{\bea}{\begin{eqnarray}}
\newcommand{\eea}{\end{eqnarray}}
\newcommand{\ba}{\begin{array}}
\newcommand{\ea}{\end{array}}

\newcommand{\bc}{\begin{center}}
\newcommand{\ec}{\end{center}}
\newcommand{\ben}{\begin{enumerate}}
\newcommand{\een}{\end{enumerate}}
\newcommand{\bfi}{\begin{figure}}
\newcommand{\efi}{\end{figure}}

\newcommand{\bq}{\begin{quote}}
\newcommand{\eq}{\end{quote}}
\newcommand{\bqu}{\begin{quotation}}
\newcommand{\equ}{\end{quotation}}
\newenvironment{emphit}{\begin{itemize}}{\end{itemize}}
\newcommand{\bemp}{\begin{emphit}}
\newcommand{\eemp}{\end{emphit}}

\newcommand{\bt}{\begin{tabular}}
\newcommand{\et}{\end{tabular}}

\newtheorem{myth}{Theorem}[section]

\newtheorem{mycor}{Corollary}[section]

\newtheorem{mydef}{Definition}[section]

\newtheorem{myrem}{Remark}[section]

\begin{document}
\date{}
\title{On Some Algebraic Properties of Generalized Groups\footnote{2000
Mathematics Subject Classification. Primary 20N99}
\thanks{{\bf Keywords and Phrases :} generalized groups}}
\author{J. O. Ad\'en\'iran, J. T. Akinmoyewa, \\
Department of Mathematics,\\
University of Agriculture, \\
Abeokuta 110101, Nigeria.\\
ekenedilichineke@yahoo.com\\
adeniranoj@unaab.edu.ng \and
A. R. T. \d S\`{o}l\'{a}r\`{i}n \\
National Mathematical Centre,\\
Federal Capital Territory,\\
P.M.B 118, Abuja, Nigeria. \\
asolarin2002@yahoo.com
\and
 T. G. Jaiy\'e\d ol\'a \footnote{corresponding author}\\
Department of Mathematics,\\
Obafemi Awolowo University,\\
Ile Ife 220005, Nigeria.\\
jaiyeolatemitope@yahoo.com\\tjayeola@oauife.edu.ng}
\maketitle
\begin{abstract}
Some results that are true in classical groups are investigated in
generalized groups and are shown to be either generally true in
generalized groups or true in some special types of generalized
groups. Also, it is shown that a Bol groupoid and a Bol quasigroup
can be constructed using a non-abelian generalized group.
\end{abstract}

\section{Introduction}
\paragraph{}
Generalized group is an algebraic structure which has a deep
physical background in the unified guage theory and has direct
relation with isotopies. Mathematicians and Physicists have been
trying to construct a suitable unified theory for twistor theory,
isotopies theory and so on. It was known that generalized groups are
tools for constructions in unified geometric theory and electroweak
theory. Electorweak theories are essentially structured on
Minkowskian axioms and gravitational theories are constructed on
Riemannian axioms. According to  Araujo et. al. \cite{gg23}, generalized group is equivalent to the notion of completely simple semigroup.

Some of the structures and properties of generalized groups have
been studied by Vagner \cite{gg25}, Molaei \cite{gg8}, \cite{gg14}, Mehrabi, Molaei and
Oloomi \cite{gg9}, Molaei and Hoseini \cite{gg15}
  and Agboola
\cite{gg10}. Smooth generalized groups were introduced in Agboola
\cite{gg11} and later on, Agboola \cite{gg12} also presented smooth
generalized subgroups while Molaei \cite{gg18} and Molaei and Tahmoresi \cite{gg24} considered the notion
of topological generalized groups. Solarin and Sharma \cite{gg16}
were able to construct a Bol loop using a group with a non-abelian
subgroup and recently, Chein and Goodaire \cite{gg22} gave a new
construction of Bol loops for odd case. Kuku \cite{gg1}, White
\cite{gg17} and Jacobson \cite{gg19} contain most of the results on
classical groups while for more on loops and their properties,
readers should check \cite{gg15,gg2,gg3,gg4,gg21,smabook1,gg20}. The aim of this study is to investigate
if some results that are true in classical group theory are also
true in generalized groups and to find a way of constructing a Bol
structure(i.e Bol loop or Bol quasigroup or Bol groupoid) using a
non-abelian generalized group.

It is shown that in a generalized group $G$,
$(a^{-1})^{-1}=a$ for all $a\in G $. In a normal generalized group
$G$, it is shown that the anti-automorphic inverse property
$(ab)^{-1}=b^{-1}a^{-1}$ for all $a,b\in G$ holds under a necessary
condition. A necessary and sufficient condition for a generalized
group(which obeys the cancellation law and in which
$e(a)=e(ab^{-1})$ if and only if $ab^{-1}=a$) to be idempotent is
established. The basic theorem used in classical groups to define
the subgroup of a group is shown to be true for generalized
groups. The kernel of any homomorphism(at a fixed point) mapping a
generalized group to another generalized group is shown to be a
normal subgroup. Furthermore, the homomorphism is found to be an
injection if and only if its kernel is the set of the identity
element at the fixed point. Given a generalized group $G$ with a
generalized subgroup $H$, it is shown that the factor set $G/H$ is
a generalized group. The direct product of two generalized group
is shown to be a generalized group. Furthermore, necessary
conditions for a generalized group $G$ to be isomorphic to the
direct product of any two abelian generalized subgroups is shown.
It is shown that a Bol groupoid can be constructed using a
non-abelian generalized group with an abelian generalized
subgroup. Furthermore, if is established that if the non-abelian
generalized group obeys the cancellation law, then a Bol
quasigroup with a left identity element can be constructed.

\section{Preliminaries}
\begin{mydef}
A generalized group $G$ is a non-empty set admitting a binary operation
called multiplication subject to the set of rules given below.
\begin{description}
\item[(i)] $(xy)z=x(yz)$ for all $x,y,z\in G$. \item[(ii)] For
each $x\in G$ there exists a unique $e(x)\in G$ such that
$xe(x)=e(x)x=x$ (existence and uniqueness of identity element).
\item[(iii)] For each $x\in G$, there exists $x^{-1}\in G$ such
that $xx^{-1}=x^{-1}x=e(x)$ (existence of inverse element).
\end{description}
\end{mydef}

\begin{mydef}
Let $L$ be a non-empty set. Define a binary operation ($\cdot $)
on $L$. If $x\cdot y\in L$ for all $x, y\in L$, $(L, \cdot )$ is
called a groupoid.

If the equations $a\cdot x=b$ and $y\cdot a=b$ have
unique solutions relative to $x$ and $y$ respectively, then $(L, \cdot )$
is called a quasigroup. Furthermore, if there exists a element
$e\in L$ called the identity element such that for all $x\in L$,
$x\cdot e=e\cdot x=x$, $(L, \cdot )$ is called a loop.
\end{mydef}

\begin{mydef}
A loop is called a Bol loop if and only if it obeys the identity
\begin{displaymath}
((xy)z)y=x((yz)y).
\end{displaymath}
\end{mydef}

\begin{myrem}
One of the most studied type of loop is the Bol loop.
\end{myrem}

\subsection{Properties of Generalized Groups}
\paragraph{}
A generalized group $G$ exhibits the following properties:
\begin{description}
\item[(i)] for each $x\in G$, there exists a unique $x^{-1}\in G$.
\item[(ii)] $e(e(x))=e(x)$ and $e(x^{-1})=e(x)$ where $x\in G$.
Then, $e(x)$ is a unique identity element of $x\in G$.
\end{description}

\begin{mydef}\label{1:4}
If $e(xy)=e(x)e(y)$ for all $x,y\in G$, then $G$ is called normal
generalized group.
\end{mydef}

\begin{myth}\label{1:5}
For each element $x$ in a generalized group $G$, there exists a
unique $x^{-1}\in G$.
\end{myth}

The next theorem shows that an abelian generalized group is a
group.

\begin{myth}\label{1:6}
Let $G$ be a generalized group and $xy=yx$ for all $x,y\in G$.
Then $G$ is a group.
\end{myth}

\begin{myth}
A non-empty subset $H$ of a generalized group $G$ is a generalized
subgroup of $G$ if and only if for all $a,b\in H$, $ab^{-1}\in H$.
\end{myth}

\paragraph{}
If $G$ and $H$ are two generalized groups and $f~:~G\to H$ is a
mapping then Mehrabi, Molaei and Oloomi \cite{gg9} called $f$ a
homomorphism if $f(ab)=f(a)f(b)$ for all $a,b\in G$.

They also stated the following results on homomorphisms of
generalized groups. These results are established in this work.

\begin{myth}\label{1:10}
Let $f~:~G\to H$ be a homomorphism where $G$ and $H$ are two
distinct generalized groups. Then:
\begin{description}
\item[(i)] $f(e(a))=e(f(a))$ is an identity element in $H$ for all
$a\in G$. \item[(ii)] $f(a^{-1})=(f(a))^{-1}$. \item[(iii)] If $K$
is a generalized subgroup of $G$, then $f(K)$ is a generalized
subgroup of $H$. \item[(iv)] If $G$ is a normal generalized group,
then the set
\begin{displaymath}
\{(e(g),f(g))~:~g\in G\}
\end{displaymath}
with the product
\begin{displaymath}
(e(a),f(a))(e(b),f(b)):=(e(ab),f(ab))
\end{displaymath}
is a generalized group denoted by $\cup f(G)$.
\end{description}
\end{myth}

\section{{\Large Main Results}}
\subsection{Results on Generalized Groups and Homomorphisms}
\begin{myth}\label{1:17}
Let $G$ be a generalized group. For all $a\in G$,
$(a^{-1})^{-1}=a$.
\end{myth}
{\bf Proof}\\
$(a^{-1})^{-1}a^{-1}=e(a^{-1})=e(a)$. Post multiplying by $a$, we
obtain
\begin{equation}\label{eq:3}
[(a^{-1})^{-1}a^{-1}]a=e(a)a.
\end{equation}
From the L. H. S.,
$(a^{-1})^{-1}(a^{-1}a)=(a^{-1})^{-1}e(a)=(a^{-1})^{-1}e(a^{-1})=(a^{-1})^{-1}e((a^{-1})^{-1})$
\begin{equation}\label{eq:4}
=(a^{-1})^{-1}.
\end{equation}
Hence from (\ref{eq:3}) and (\ref{eq:4}), $(a^{-1})^{-1}=a$.

\begin{myth}\label{1:18}
Let $G$ be a generalized group in which the left cancellation law
holds and $e(a)=e(ab^{-1})$ if and only if $ab^{-1}=a$. $G$ is a
idempotent generalized group if and only if
$e(a)b^{-1}=b^{-1}e(a)~\forall~a,b\in G$.
\end{myth}
{\bf Proof}\\
$e(a)b^{-1}=b^{-1}e(a)\Leftrightarrow
(ae(a))b^{-1}=ab^{-1}e(a)\Leftrightarrow
ab^{-1}=ab^{-1}e(a)\Leftrightarrow e(a)=e(ab^{-1})\Leftrightarrow
ab^{-1}=a\Leftrightarrow ab^{-1}b=ab\Leftrightarrow
ae(b)=ab\Leftrightarrow a^{-1}ae(b)=a^{-1}ab\Leftrightarrow
e(a)e(b)=e(a)b\Leftrightarrow e(b)=b\Leftrightarrow b=bb$.

\begin{myth}\label{1:19}
Let $G$ be a normal generalized group in which
$e(a)b^{-1}=b^{-1}e(a)~\forall~a,b\in G$. Then,
$(ab)^{-1}=b^{-1}a^{-1}~\forall~a,b\in G$.
\end{myth}
{\bf Proof}\\
Since $(ab)^{-1}(ab)=e(ab)$, then by multiplying both sides of the
equation on the right by $b^{-1}a^{-1}$ we obtain
\begin{equation}\label{eq:1}
[(ab)^{-1}ab]b^{-1}a^{-1}=e(ab)b^{-1}a^{-1}.
\end{equation}
So,
$[(ab)^{-1}ab]b^{-1}a^{-1}=(ab)^{-1}a(bb^{-1})a^{-1}=(ab)^{-1}a(e(b)a^{-1})=(ab)^{-1}(aa^{-1})e(b)=(ab)^{-1}(e(a)e(b))
=(ab)^{-1}e(ab)=(ab)^{-1}e((ab)^{-1})$
\begin{equation}\label{eq:2}
=(ab)^{-1}.
\end{equation}
Using (\ref{eq:1}) and (\ref{eq:2}), we get
$[(ab)^{-1}ab]b^{-1}a^{-1}=(ab)^{-1}\Rightarrow
e(ab)(b^{-1}a^{-1})=(ab)^{-1}\Rightarrow (ab)^{-1}=b^{-1}a^{-1}$.

\begin{myth}\label{1:20}
Let $H$ be a non-empty subset of a generalized group $G$. The
following are equivalent.
\begin{description}
\item[(i)] $H$ is a generalized subgroup of $G$. \item[(ii)] For
$a,b\in H$, $ab^{-1}\in H$. \item[(iii)] For $a,b\in H$, $ab\in H$
and for any $a\in H$, $a^{-1}\in H$.
\end{description}
\end{myth}
{\bf Proof}\\
\paragraph{(i)$\Rightarrow$ (ii)}
If $H$ is a generalized subgroup of $G$ and $b\in G$, then
$b^{-1}\in H$. So by closure property, $ab^{-1}\in H~\forall~a\in
H$.
\paragraph{(ii)$\Rightarrow$ (iii)}
If $H\ne \phi$, and $a,b\in H$, then we have $bb^{-1}=e(b)\in H$,
$e(b)b^{-1}=b^{-1}\in H$ and $ab=a(b^{-1})^{-1}\in H$ i.e $ab\in
H$.
\paragraph{(iii)$\Rightarrow$ (i)}
$H\subseteq G$ so $H$ is associative since $G$ is associative.
Obviously, for any $a\in H$, $a^{-1}\in H$. Let $a\in H$, then
$a^{-1}\in H$. So, $aa^{-1}=a^{-1}a=e(a)\in H$. Thus, $H$ is a
generalized subgroup of $G$.

\begin{myth}\label{1:21}
Let $a\in G$ and $f~:~G\to H$ be an homomorphism. If $\ker f$ at
$a$ is denoted by
\begin{displaymath}
\ker f_a=\{x\in G~:~f(x)=f(e(a))\}.
\end{displaymath}
Then,
\begin{description}
\item[(i)] $\ker f_a\triangleleft G$. \item[(ii)] $f$ is a
monomorphism if and only if $\ker f_a=\{e(a)~:~\forall~a\in G\}$.
\end{description}
\end{myth}
{\bf Proof}\\
\begin{description}
\item[(i)] It is necessary to show that $\ker f_a\le G$. Let
$x,y\in \ker f_a\le G$, then
$f(xy^{-1})=f(x)f(y^{-1})=f(e(a))(f(e(a)))^{-1}=f(e(a))f(e(a)^{-1})=f(e(a))f(e(a))=f(e(a))$.
So, $xy^{-1}\in \ker f_a$. Thus, $\ker f_a\le G$. To show that
$\ker f_a\triangleleft G$, since $y\in \ker f_a$, then by the
definition of $\ker f_a$,
$f(xyx^{-1})=f(x)f(y)f(x^{-1})=f(e(a))f(e(a))f(e(a))^{-1}=f(e(a))f(e(a))f(e(a))=f(e(a))\Rightarrow
xyx^{-1} \ker f_a$. So, $\ker f_a\triangleleft G$. \item[(ii)]
$f~:~G\to H$. Let $\ker f_a=\{e(a)~:~\forall~a\in G\}$ and
$f(x)=f(y)$, this implies that
$f(x)f(y)^{-1}=f(y)f(y)^{-1}\Rightarrow
f(xy^{-1})=e(f(y))=f(e(y))\Rightarrow xy^{-1}\in
\ker f_y\Rightarrow$
\begin{equation}\label{eq:5}
xy^{-1}=e(y)
\end{equation}
and $f(x)f(y)^{-1}=f(x)f(x)^{-1}\Rightarrow
f(xy^{-1})=e(f(x))=f(e(x))\Rightarrow xy^{-1}\in
\ker f_x\Rightarrow$
\begin{equation}\label{eq:6}
xy^{-1}=e(x).
\end{equation}
Using (\ref{eq:5}) and (\ref{eq:6}),
$xy^{-1}=e(y)=e(x)\Leftrightarrow x=y$. So, $f$ is a
monomorphism.\\

Conversely, if $f$ is mono, then $f(y)=f(x)\Rightarrow y=x$. Let
$k\in \ker f_a~\forall~a\in G$. Then, $f(k)=f(e(a))\Rightarrow
k=e(a)$. So, $\ker f_a=\{e(a)~:~\forall~a\in G\}$.
\end{description}

\begin{myth}\label{1:21.2}
Let $G$ be a generalized group and $H$ a generalized subgroup of
$G$. Then $G/H$ is a generalized group called the quotient or
factor generalized group of $G$ by $H$.
\end{myth}
{\bf Proof}\\
It is necessary to check the axioms of generalized group on $G/H$.
\begin{itemize}
\item[Associativity] Let $a,b,c\in G$ and $aH,bH,cH\in G/H$. Then
$aH(bH\cdot cH)=(aH\cdot bH)cH$, so associativity law holds.
\item[Identity] If $e(a)$ is the identity element for each $a\in
G$, then $e(a)H$ is the identity element of $aH$ in $G/H$ since
$e(a)H\cdot aH=e(a)\cdot aH=aH\cdot e(a)=aH$. Therefore identity
element exists and is unique for each elements $aH$ in $G/H$.
\item[Inverse]
$(aH)(a^{-1}H)=(aa^{-1})H=e(a)H=(a^{-1}a)H=(a^{-1}H)(aH)$ shows
that $a^{-1}H$ is the inverse of $aH$ in $G/H$.
\end{itemize}
So the axioms of generalized group are satisfied in $G/H$.

\begin{myth}\label{1:21.1}
Let $G$ and $H$ be two generalized groups. The direct product of
$G$ and $H$ denoted by
\begin{displaymath}
G\times H=\{(g,h)~:~g\in G~\textrm{and}~h\in H\}
\end{displaymath}
is a generalized group under the binary operation $\circ$ such
that
\begin{displaymath}
(g_1,h_1)\circ (g_2,h_2)=(g_1g_2,h_1h_2).
\end{displaymath}
\end{myth}
{\bf Proof}\\
This is achieved by investigating the axioms of generalized group
for the pair $(G\times H, \circ )$.

\begin{myth}\label{1:22}
Let $G$ be a generalized group with two abelian generalized
subgroups $N$ and $H$ of $G$ such $G=NH$. If $N\subseteq COM(H)$
or $H\subseteq COM(N)$ where $COM(N)$ and $COM(H)$ represent the
commutators of $N$ and $H$ respectively, then $G\cong N\times H$.
\end{myth}
{\bf Proof}\\
Let $a\in G$. Then $a=nh$ for some $n\in N$ and $h\in H$. Also,
let $a=n_1h_1$ for some $n_1\in N$ and $h_1\in H$. Then
$nh=n_1h_1$ so that $e(nh)=e(n_1h_1)$, therefore $n=n_1$ and
$h=h_1$. So that $a=nh$ is unique.

Define $f~:~G\to H$ by $f(a)=(n,h)$ where $a=nh$. This function is
well defined in the previous paragraph which also shows that $f$
is a one-one correspondence. It remains to check that $f$ is a group
homomorphism.

Suppose that $a=nh$ and $b=n_1h_1$, then $ab=nhn_1h_1$ and
$hn_1=n_1h$. Therefore,
$f(ab)=f(nhn_1h_1)=f(nn_1hh_1)=(nn_1,hh_1)=(n,h)(n_1,h_1)=f(a)f(b)$.
So, $f$ is a group homomorphism. Hence a group isomorphism since
it is a bijection.

\subsection{Construction of Bol Algebraic Structures}

\begin{myth}\label{1:24}
Let $H$ be a subgroup of a non-abelian generalized group $G$ and
let $A=H\times G$. For $(h_1,g_1),(h_2,g_2)\in A$, define
\begin{displaymath}
(h_1,g_1)\circ (h_2,g_2)=(h_1h_2,h_2g_1h_2^{-1}g_2)
\end{displaymath}
then $(A,\circ )$ is a Bol groupoid.
\end{myth}
{\bf Proof}\\
Let $x,y,z\in A$. By checking, it is true that $x\circ (y\circ
z)\ne (x\circ y)\circ z$. So, $(A,\circ )$ is non-associative. $H$
is a quasigroup and a loop(groups are quasigroups and loops) but
$G$ is neither a quasigroup nor a loop(generalized groups are
neither quasigroups nor a loops) so $A$ is neither a quasigroup
nor a loop
but is a groupoid because $H$ and $G$ are groupoids.\\

Let us now verify the Bol identity:
\begin{displaymath}
((x\circ y)\circ z)\circ y=x\circ ((y\circ z)\circ y)
\end{displaymath}
\begin{displaymath}
\textrm{L. H. S.}\qquad =((x\circ y)\circ z)\circ
y=(h_1h_2h_3h_2,h_2h_3h_2g_1h_2^{-1}g_2h_3^{-1}g_3h_2^{-1}g_2).
\end{displaymath}
\begin{displaymath}
\textrm{R. H. S.}\qquad = x\circ ((y\circ z)\circ
y)=(h_1h_2h_3h_2,h_2h_3h_2g_1h_2^{-1}(h_3^{-1}h_2^{-1}h_2h_3)g_2h_3^{-1}g_3h_2^{-1}g_2)=
\end{displaymath}
\begin{displaymath}
(h_1h_2h_3h_2,h_2h_3h_2g_1h_2^{-1}g_2h_3^{-1}g_3h_2^{-1}g_2).
\end{displaymath}
So, L. H. S.=R. H. S.. Hence, $(A,\circ )$ is a Bol groupoid.

\begin{mycor}\label{1:25}
Let $H$ be a abelian generalized subgroup of a non-abelian
generalized group $G$ and let $A=H\times G$. For
$(h_1,g_1),(h_2,g_2)\in A$, define
\begin{displaymath}
(h_1,g_1)\circ (h_2,g_2)=(h_1h_2,h_2g_1h_2^{-1}g_2)
\end{displaymath}
then $(A,\circ )$ is a Bol groupoid.
\end{mycor}
{\bf Proof}\\
By Theorem~\ref{1:6}, an abelian generalized group is a group, so
$H$ is a group. The rest of the claim follows from
Theorem~\ref{1:24}.

\begin{mycor}\label{1:26}
Let $H$ be a subgroup of a non-abelian generalized group $G$ such
that $G$ has the cancellation law and let $A=H\times G$. For
$(h_1,g_1),(h_2,g_2)\in A$, define
\begin{displaymath}
(h_1,g_1)\circ (h_2,g_2)=(h_1h_2,h_2g_1h_2^{-1}g_2)
\end{displaymath}
then $(A,\circ )$ is a Bol quasigroup with a left identity
element.
\end{mycor}
{\bf Proof}\\
The proof of this goes in line with Theorem~\ref{1:24}. A groupoid
which has the cancellation law is a quasigroup, so $G$ is
quasigroup hence $A$ is a quasigroup. Thus, $(A,\circ )$ is a Bol
quasigroup with a left identity element since by kunen
\cite{kun1}, every quasigroup satisfying the right Bol identity
has a left identity.

\begin{mycor}\label{1:27}
Let $H$ be a abelian generalized subgroup of a non-abelian
generalized group $G$ such that $G$ has the cancellation law and
let $A=H\times G$. For $(h_1,g_1),(h_2,g_2)\in A$, define
\begin{displaymath}
(h_1,g_1)\circ (h_2,g_2)=(h_1h_2,h_2g_1h_2^{-1}g_2)
\end{displaymath}
then $(A,\circ )$ is a Bol quasigroup with a left identity
element.
\end{mycor}
{\bf Proof}\\
By Theorem~\ref{1:6}, an abelian generalized group is a group, so
$H$ is a group. The rest of the claim follows from
Theorem~\ref{1:26}.

\end{document}